\documentclass[openacc]{rstransa}

\usepackage{commath}

\usepackage{tikz}
\usetikzlibrary{fit}
\newcommand\tm[2][]{\tikz[overlay,remember picture,baseline=(#1.base),inner sep=0pt]\node(#1){$#2$};}
\usepackage{accents}

\begin{document}

\title{On Incorporating Forecasts into Linear State Space Model Markov Decision Processes}

\author{
Jacques A. de Chalendar$^{1}$ and Peter W. Glynn$^{2}$}

\address{$^{1}$Department of Energy Resources Engineering, Stanford University, Stanford CA 94305-2205, USA.\\
$^{2}$Department of Management Sciences \& Engineering, Stanford University, Stanford CA 94305-2205, USA.}

\subject{xxxxx, xxxxx, xxxx}

\keywords{Markov Decision Processes, weather forecasts, forecast models, stochastic control of energy systems}

\corres{Jacques de Chalendar\\
\email{jdechalendar@stanford.edu}}

\begin{abstract}
Weather forecast information will very likely find increasing application in the control of future energy systems. In this paper, we introduce an augmented state space model formulation with linear dynamics, within which one can incorporate forecast information that is dynamically revealed alongside the evolution of the underlying state variable. We use the martingale model for forecast evolution (MMFE) to enforce the necessary consistency properties that must govern the joint evolution of forecasts with the underlying state. The formulation also generates jointly Markovian dynamics that give rise to Markov decision processes (MDPs) that remain computationally tractable. This paper is the first to enforce MMFE consistency requirements within an MDP formulation that preserves tractability.
\end{abstract}


\begin{fmtext}

\end{fmtext}


\maketitle
\section{Introduction}
Forecasts are ubiquitous in energy system control problems and there is reason to believe their importance will only grow, \textit{e.g.} in the fast changing electric sector. This is especially true for forecasts that provide weather-related information, as weather patterns have a strong impact on energy demand and increasingly on (renewable) energy production. The meteorological community has made significant progress in that field over the past decades and can now offer several advantages over purely statistical models \cite{bauer2015quiet}. In a recent review of forecasting for renewable energy \cite{sweeney2020future}, the authors note a rise in demand for probabilistic forecasting and that the typical renewable energy production use cases for weather forecasts correspond to timescales of weeks to years, or hours to a day ahead.
\par Much deeper transformations are in store once very high penetrations of renewables are reached. A growing number of applications will require storing electricity with durations from 10 to 100 hours \cite{shaner2018geophysical, ziegler2019storage, albertus2020long}. Wind generation is one example for which low availability levels can be observed for several consecutive days. The need for longer duration storage will appear more clearly with penetrations of > 70\% wind and solar generation on a regional grid – see e.g. Figure 1 from \cite{albertus2020long}. As soon as longer duration storage becomes available, we will correspondingly need new energy management strategies. At these timescales, conditioning on forecasts, in particular for weather variables, will very likely have large impacts on problems of decision making under uncertainty.
\par Using forecast information in the context of control problems is a difficult general problem that implicitly appears in many real-life applications. In sequential decision problems, it is often the case that exogenous forecast information is presented to the controller at regular intervals. Given the key role of Markov decision processes (MDPs) in the computation of optimal policies in such settings, a full accounting of the impact of future forecast information requires introducing the forecasts into the Markov state variable, thereby leading to potential high-dimensional state representations. Another fundamental issue relates to the fact that the forecasts should be ``compatible'' with the state variable that is being forecast (\textit{e.g.} weather), so that the forecasted state variable (\textit{i.e.} the state variable for which forecasts are available) and the forecasts themselves should exhibit self-consistent dynamics. To gain some appreciation for this issue, note that the $s$-period forecast must contain information that implicitly ``peeks'' $s$ periods into the future of the underlying state space model, so that the $s$-period forecast implicitly constrains the dynamics of the underlying model over the next $s$ periods. These constraints need to be built into the joint dynamics in such a way that the Markov structure is preserved. The preservation of Markov structure is critical if we wish to be able to compute optimal policies via the use of MDP-based theory and algorithms.
\par In this paper, we utilize the martingale model for forecast evolution (MMFE) as a vehicle for imposing the appropriate mathematical consistency between the dynamics of the forecasts and the forecasted state variables. The MMFE framework was introduced and developed by \cite{hausman1969sequential, graves1986two, heath1994modeling} and has since been utilized extensively by the inventory control and supply chain management community (see \textit{e.g} \cite{altug2011inventory} and the references therein). Applications of the MMFE and studies on the impact of forecasts on decision making can also be found in the energy community, \textit{e.g.} for hydraulic reservoir management \cite{zhao2013generalized} or wind energy integration \cite{bitar2012bringing, nair2014energy}.
\par To our knowledge, this paper is the first to rigorously introduce forecast model consistency into MDPs, specifically in the context of linear state space models. This work gives us the first principled and mathematically consistent framework for the incorporation of forecasts into MDPs in the setting of state space models with linear dynamics, and uses no ad hoc elements to add forecast information into the MDP setting. Linear state space models are widely applied across many disciplines, and can even represent the linearized dynamics associated with nonlinear structure \cite{callier2012linear}.
\par We note that forecasts must depend on a richer information filtration than that associated with the forecasted state variable, because an optimal MDP policy computed from the forecasted state variable already fully utilizes all the information associated with the forecasted state variable's filtration. In our setting, the extra information that enters the forecasts is the meteorological data available to the forecasters that is unobserved by the energy system manager. Thus, a key contribution of our paper is the development of an MMFE framework in which one can rigorously discuss, via the use of the language of $\sigma$-algebras, the different ways in which forecast information can be incorporated into an MDP framework. In our carefully chosen formulation, each of these different approaches for incorporating MMFE forecasts leads to a different, but computationally tractable, MDP.
\par Our first new MDP (Section~\ref{sec:4}) incorporates the ``static'' forecast information that is available to the decision-maker at the beginning of the decision horizon, and leads to an MDP that has the same state space as for the forecasted state variable, but with transition probabilities that are non-stationary as a consequence of the initial set of forecasts. The model has the property that when one conditions the future dynamics of the forecasted state variable on the forecast information available at the beginning of the decision horizon, the Markov structure is preserved with no need to increase the dimensionality of the state representation. In Section~\ref{sec:5}, we develop an MDP in which the forecasts are dynamically updated over time, along with the forecasted state variable. Thus, this formulation explicitly models the additional forecasting information that is revealed to the decision-maker over time. In this dynamic forecasting formulation, one needs to expand the state space of the MDP to incorporate the forecast evolution, but the MDP has stationary transition probabilities. Our final MDP is a formulation in which new $r$-period lookahead forecasts are made available to the decision-maker over time, in addition to an extended set of static forecasts that provides forecast information more than $r$ periods into the future. This MDP that combines both static and dynamic forecasts is introduced in Section~\ref{sec:6}, and leads to both an enlarged state space and non-stationary transition probabilities.
\par An alternative means of utilizing the availability of forecasts in the control setting is to apply the ideas of model predictive control (MPC). MPC has become a standard tool for many industrial applications and provides a practical way of dealing with forecasts \cite{qin2003survey, mayne2014model, mesbah2016stochastic}. In this approach, one uses the forecasts available to solve a sequence of MDP formulations over time. At each decision epoch, a conventional MDP that incorporates the forecasts available is solved, the optimal first period action is taken, and this process is repeated at the next decision epoch. In particular, the MDPs that are used by MPC at each decision epoch do not explicitly model, within the MDP, the fact that the decision-maker will have available a new set of forecasts at each future decision epoch within the decision horizon associated with the MDP. This is also the case in adaptations of standard MPC to the setting where the forecasts contain probabilistic information \cite{mesbah2016stochastic}. In contrast, the MDPs introduced by this paper model the fact that the forecasts are continually ``refreshed'' over the MDP's decision horizon, and do so via a formulation that preserves the computational tractability of the MDP.
\par In Section~\ref{sec:7}, we introduce a simple energy system control model that controls interior building temperatures in an external weather environment for which forecasts are available. In the presence of a quadratic cost structure, we are able to use the existing linear-quadratic stochastic control theory to compute the optimal value associated with MDPs for the energy system in which no forecast information is available and also the optimal value for the dynamic forecasting MDP of Section~\ref{sec:5}. This allows us to analyze the degree of improvement that can be obtained by incorporating dynamic forecast information into the MDP formulation in the context of our simple energy system example. Section~\ref{sec:8} concludes the paper with a discussion of additional research questions that this paper motivates.

\section{MDP's with no forecasts}
\label{sec:2}
In this section, we review the basic MDP framework that can be used when making sequential decisions involving an energy system that is affected by the weather. Our formulation here does not take advantage of any forecast information that may be available. We model the dynamics in discrete time, and take the view that the weather variables $W_n$ at time $n$ can be represented by an $\mathbb{R}^d$-valued random variable (rv). Given our MDP modeling perspective, we assume that $(W_n:n \in \mathbb{Z})$ is a stationary $\mathbb{R}^d$-valued stochastic process that enjoys the Markov property, so that
\begin{equation}
    W_{n+1} = f(W_n, Z_{n+1})
    \label{eq:2.1}
\end{equation}
for $n\in\mathbb{Z}$, where $(Z_n:n \in \mathbb{Z})$ is an independent and identically distributed (iid) sequence of $\mathbb{R}^{m_1}$-valued rv's and $f$ is a (deterministic) mapping from $\mathbb{R}^d\times\mathbb{R}^{m_1}$ into $\mathbb{R}^d$. The stationarity is intended here to simplify the exposition and is (at best) only approximately valid in the weather setting. For example, in examining daily weather records, it may be that such time series look approximately stationary over time scales of (say) one month. Given that our decision horizon is typically much shorter than a month, the stationarity assumption will often be a reasonable one in practice.
\par For the energy system's control, we model its state evolution via an $\mathbb{R}^l$-valued sequence $(X_n:n\in\mathbb{Z})$ for which
\begin{equation}
    X_{n+1} = \phi(X_n, A_n, W_{n+1}, V_{n+1})
    \label{eq:2.2}
\end{equation}
for $n\in\mathbb{Z}$, where $(V_n:n\in\mathbb{Z})$ is an iid sequence of $\mathbb{R}^{m_2}$-valued rv's independent of the $Z_j$'s, $A_n$ is an $\mathcal{A}$-valued action taken at time $n$, and $\phi$ is a deterministic mapping. The action $A_n$ must be adapted to the history $\mathcal{F}_n\triangleq\sigma((W_j,X_j):j\leq n)$, so that it can depend only on previously observed values of the weather and control system state. The joint dynamics \ref{eq:2.1} and \ref{eq:2.2} assume (reasonably) that the weather affects the control system dynamics, but not vice-versa.
\par We now describe the dynamic program (DP) backwards dynamic recursion that is commonly used to compute the optimal $A^*_j$'s, when optimizing the control of such an energy system over a finite horizon $[0, t+1)$. Throughout this paper, we take $n=0$ as the time at which the sequence of control actions will be computed. Suppose that our goal is to minimize the total expected cost of running the energy system over $[0, t+1)$, namely
\begin{equation}
    \mathbb{E}\left[\sum_{j=0}^tc(X_j,A_j,W_{j+1})|X_0=x,W_0=w\right]
\end{equation}
over all $F_{j}$-adapted controls $(A_j:0\leq j\leq t)$. Here $c(X_j, A_j, W_{j+1})$ represents the one-period cost for running the energy system over $[j, j+1)$. For an (appropriately integrable) function $h$ with domain $\mathbb{R}^l\times\mathbb{R}^d$, define the operator
\begin{equation}
    (P_ah)(x,w)=\int_{\mathbb{R}^{m_1}}\int_{\mathbb{R}^{m_2}}h(\phi(x, a, f(w,z),v),f(w,z))\mathbb{P}(Z_1 \in dz)\mathbb{P}(V_1 \in dv)
\end{equation}
and set
\begin{equation}
    \tilde{c}(x, a, w) = \int_{\mathbb{R}^{m_1}} c(x, a, f(w, z))\mathbb{P}(Z_1 \in dz)
\end{equation}
The DP value functions $(v_i(\cdot):0\leq i\leq t)$ are then computed via the recursion
\begin{equation}
    v_i(x,w) = \min_{a\in\mathcal{A}}\left[\tilde{c}(x, a, w) + (P_av_{i+1})(x, w)\right]
    \label{eq:2.6}
\end{equation}
for $0\leq i< t$, subject to the terminal condition
\begin{equation}
    v_t(x,w) = \min_{a\in\mathcal{A}} \tilde{c}(x, a, w).
    \label{eq:2.7}
\end{equation}
Assuming that $v_t, v_{t-1}, \cdots,v_0$ are recursively computed via \ref{eq:2.6} and \ref{eq:2.7} , we then select $a^*_i(x, w)$ ($a^*_t(x, w)$) as any minimizer (assumed to exist) of the right-hand side of \ref{eq:2.6} (\ref{eq:2.7}), and put $A_i^*=a_i^*(X_i, W_i)$ for $0\leq i \leq t$. Under suitable integrability hypotheses, it is well known that $(A_i^*: 0\leq i \leq t)$ is then the desired cost-minimizing adapted optimal control, see \textit{e.g.} \cite{bertsekas2004stochastic}.

\section{The mathematical structure of forecasts}
\label{sec:3}
In order to build forecast information into the Markov model of Section~\ref{sec:2}, we review the mathematical structure of forecasts, so that we can ensure that the model combining Markovian state dynamics and forecast information respects the appropriate mathematical constraints. To this end, we assume that $\mathbb{E}||W_n||^2< \infty$ (where $||\cdot||$ is the Euclidian norm). We model the (point) \textit{forecast} $F_{n|k}$ of $W_n$ available at time $k\leq n$ as the rv $\mathbb{E}[W_n|\mathcal{G}_k]$, where $\mathcal{G}_k$ is a $\sigma$-algebra representing the information available to the forecaster. Since weather forecasters have available vastly more weather information than does the energy system manager, we expect that $\mathcal{G}_k$ represents a strictly richer ``information set'' than $\mathcal{F}^W_k\triangleq\sigma(W_j:j\leq k)$. Consequently, we require $\mathcal{G}_k$ to be strictly larger than $\mathcal{F}^W_k$. In fact, if $\mathcal{G}_k = \mathcal{F}^W_k$, the availability of forecasts will offer no advantage over the optimal control $(A^*_j: 0\leq j \leq t)$ computed in Section~\ref{sec:2}, since that policy is already guaranteed to be optimal over all $\mathcal{F}^W_k$-adapted policies.
\par We note that $F_{n|n}=W_n$ and that the tower property of conditional expectation implies that
\begin{equation}
    \begin{split}
        \mathbb{E}[F_{n|k+1}|\mathcal{G}_k] &= \mathbb{E}\left[\mathbb{E}[W_n|\mathcal{G}_{k+1}]|\mathcal{G}_k\right] \\
        &=\mathbb{E}[W_n|\mathcal{G}_{k}]\\
        &=F_{n|k},
    \end{split}
\end{equation}
so that $(F_{n|k}:k\leq n)$ is a \textit{martingale} (in $k$) adapted to the $\mathcal{G}_k$'s for each fixed $n\in\mathbb{Z}$. For $k\leq n$, let
\begin{equation}
    D_{n|k}=F_{n|k}-F_{n|k-1}
\end{equation}
be the $k$'th \textit{martingale difference} associated with the martingale $(F_{n|k}:k\leq n)$. The square integrability of the $W_n$'s implies that $D_{n|k}D_{m|j}$ is integrable and
\begin{equation}
    \mathbb{E}[D_{n|k}D_{m|j}|\mathcal{G}_j] = D_{m|j}\mathbb{E}\left[\mathbb{E}[D_{n|k}|\mathcal{G}_k]|\mathcal{G}_j\right]=0
\end{equation}
for $j\leq k$, so that
\begin{equation}
    \mathbb{E}D_{n|k}D_{m|j} = 0
\end{equation}
for $j\neq k$ and $n\geq k, m\geq j$. This \textit{orthogonality} of $D_{n|k}$ and $D_{m|j}$ is a key property of such martingale differences. As was discussed in the Introduction, the fact that such martingale structure is a reasonable requirement to impose on forecasts has been noted previously \cite{hausman1969sequential, graves1986two, heath1994modeling, nair2014energy}.

\section{MDP's incorporating a static forecast}
\label{sec:4}
We now wish to build a tractable model under which $(W_k: 0\leq k \leq t+1)$ evolves over the decision horizon, conditional on the forecasts $(F_{n|0}:n\geq 0)$ available at the outset of the decision interval.
\par To this end, let $\mathcal{K}_n \triangleq\sigma(F_{m+j|m}:j\in\mathbb{Z}_+, m\leq n)$ denote the $\sigma$-algebra associated with the forecasts collected by time $n$, and note that $\mathcal{K}_n\subseteq\mathcal{G}_n$, the $\sigma$-algebra associated with all the information observed by the forecaster by time $n$. We now wish to construct an MDP formulation appropriate to decision-making by the energy system manager when she has access to the information available both in $\mathcal{F}_n$ and $\mathcal{K}_0$. In other words, her decision at time $n$ must be $\mathcal{F}_n\vee\mathcal{K}_0$ adapted, where $\mathcal{B}_1\vee\mathcal{B}_2\vee\dots\vee\mathcal{B}_l$ is our notation for the smallest $\sigma$-algebra containing $\mathcal{B}_1, \mathcal{B}_2, \dots, \mathcal{B}_l$. We call this a \textit{static forecast} formulation, since the decision maker only uses the forecasts available at time 0 in making decisions.
\par In particular, we shall build a model under which the (conditional on $\mathcal{G}_0$) Markov property
\begin{equation}
    \mathbb{P}(W_{n+1}\in\cdot|\mathcal{G}_0, W_j:j\leq n) = \mathbb{P}(W_{n+1}\in\cdot|\mathcal{K}_0, W_n)
    \label{eq:4.1}
\end{equation}
holds for $0\leq n\leq t$. This ensures that
\begin{equation}
    \mathbb{P}(W_{n+1}\in\cdot|\mathcal{K}_0,W_j:0\leq j\leq n) = \mathbb{P}(W_{n+1}\in\cdot|\mathcal{K}_0,W_n).
\end{equation}
\par We will now formulate a flexible model that satisfies both the ordinary Markov property (as expressed through the recursion \ref{eq:2.1}) and the conditional Markov property (as expressed through \ref{eq:4.1}). In particular, we now specialize the stochastic recursion \ref{eq:2.1} to a \textit{linear state space model} of the form
\begin{equation}
    W_{n+1} = GW_n+Z_{n+1}
    \label{eq:4.2}
\end{equation}
where $G$ is a deterministic $d\times d$ matrix having spectral radius less than 1, and for which the $Z_i$'s are iid $\mathbb{R}^d$-valued rv's for which $\mathbb{E}||Z_1||^2<\infty$. Then, $\mathbb{E}W_n=(I-G)^{-1}\mathbb{E}Z_1$ for $n\in\mathbb{Z}$. Let $\tilde{W}_n=W_n-\mathbb{E}W_0$ and $\tilde{Z}_n=Z_n-\mathbb{E}Z_1$ and note that \ref{eq:4.2} implies that
\begin{equation}
    \tilde{W}_{n+1}=G\tilde{W}_n+\tilde{Z}_{n+1}
\end{equation}
for $n\in\mathbb{Z}$. We further assume that for each $n\in\mathbb{Z}$, we can write $\tilde{Z}_n$ in the form
\begin{equation}
    \tilde{Z}_n=\sum_{j=0}^\infty\epsilon_n(n-j),
    \label{eq:4.4}
\end{equation}
where the sum in \ref{eq:4.4} is assumed to converge a.s. and in mean square. The family of rv's $(\epsilon_n(j), j\leq n, n\in\mathbb{Z})$ is assumed to satisfy:
\renewcommand{\labelenumi}{A\arabic{enumi}.}
\begin{enumerate}
    \item The collection $(\epsilon_n(n-j):n\in\mathbb{Z}, j\in\mathbb{Z}_+)$ is a family of independent mean zero square integrable rv's.\\
    \item $\epsilon_n(n-j) \stackrel{\mathcal{D}}{=}\epsilon_0(-j)$, for $n\in\mathbb{Z}$ (where $\stackrel{\mathcal{D}}{=}$ means ``equality in distribution'').
\end{enumerate}
\paragraph{Remark}
The $\epsilon_n(k)$ disturbance models the information gathered by the forecaster at time $k$ that is relevant to the forecast for time $n$. In view of this interpretation, it is natural that we then ``model'' the $\sigma$-algebra $\mathcal{G}_n$ of Section~\ref{sec:3} as $\mathcal{G}_n\triangleq\sigma(\epsilon_{m}(j):j\leq n, m\geq j)$ in the context of this state space model. In this case $\mathcal{G}_n=\mathcal{K}_n$, as we will see later in this section, although we note that in general $\mathcal{G}_n$ could be strictly richer than $\mathcal{K}_n$. We further note that A2 implies that the distribution for $\epsilon_n(k)$ only depends on $n-k$. A1 and A2 ensure that $(\tilde{Z}_n:n\in\mathbb{Z})$ is an iid sequence of mean zero square integrable rv's.
\\\\
If we set $\mathcal{H}_n\triangleq\sigma(\epsilon_n(m-j): m\leq n, j\in\mathbb{Z}_+)$, we note that $\mathcal{F}_n^W\triangleq\sigma(W_j:j\leq n) \subseteq \mathcal{H}_n$ and that the independence of $(\epsilon_{n+1}(n+1-j):j\in\mathbb{Z}_+)$ from $\mathcal{H}_n$ ensures that
\begin{equation}
    \mathbb{P}(W_{n+1}\in\cdot|\mathcal{H}_n) = \mathbb{P}(W_{n+1}\in\cdot|W_n).
\end{equation}
It follows that the policy $(A_n^*: 0 \leq n \leq t)$ computed in Section~\ref{sec:2} is optimal not only over the $\mathcal{F}_n$-adapted policies but also over the $\mathcal{H}_n\vee\mathcal{F}_n$-adapted policies.
\par Furthermore, for $k\leq n$,
\begin{equation}
    \begin{split}
        W_n&=\mathbb{E}W_0 + G^{n-k}\tilde{W}_k + \sum_{i=0}^{n-k-1}G^i\tilde{Z}_{n-i}\\
        &=\mathbb{E}W_0 + G^{n-k}\tilde{W}_k + \sum_{i=0}^{n-k-1}G^i \sum_{j=0}^{\infty}\epsilon_{n-i}(n-i-j)\\
        &=\mathbb{E}W_0 + G^{n-k}\tilde{W}_k + \sum_{i=0}^{n-k-1}G^i \sum_{r=-\infty}^{n-k-1}\epsilon_{n-i}(r)\\
        &=\mathbb{E}W_0 + G^{n-k}\tilde{W}_k + \sum_{r=-\infty}^{k} \sum_{i=0}^{n-k-1} G^i \epsilon_{n-i}(r) + \sum_{r=k+1}^{n} \sum_{i=0}^{n-r} G^i \epsilon_{n-i}(r).
    \end{split}
    \label{eq:4.5}
\end{equation}
\par We recall that $\mathcal{G}_n$ contains the sequence of rv's $(\epsilon_m(j):j\leq n, m>n)$ that are independent of $\mathcal{H}_n$ (and hence $\mathcal{F}_n$). This represents the additional information available to the forecaster about the weather in future time periods that goes beyond the predictive information present in observing $W_n$ that is locally available to the energy system manager. Figures~\ref{fig:1} and~\ref{fig:2} illustrate the differences in the information sets $\mathcal{H}_n$ and $\mathcal{G}_n$.

\begin{figure}
    \centering
    \begin{equation*}
    \begin{matrix}
    \tm[a]{\ddots} &                   &                   &        &                     \\
    \cdots & \epsilon_{-1}(-1) &                   &        &                     \\
    \cdots & \epsilon_{0}(-1)  & \epsilon_0(0)     &        &                    \\
           & \vdots            & \vdots            & \ddots &                    \\
    \cdots & \epsilon_{n-1}(-1)  & \epsilon_{n-1}(0) & \cdots & \epsilon_{n-1}(n-1)\\
    \tm[b]{\cdots} & \tm[b1]{\epsilon_{n}(-1)}  & \epsilon_{n}(0) & \cdots & \epsilon_{n}(n-1)& \tm[c]{\epsilon_{n}(n)}\\
    \cdots & \epsilon_{n+1}(-1)  & \epsilon_{n+1}(0) & \cdots & \epsilon_{n+1}(n-1)& \epsilon_{n+1}(n) & \epsilon_{n+1}(n+1)\\
     & \vdots  & \vdots & & \vdots& \vdots & \vdots
    \end{matrix}    
    \end{equation*}
    
    \begin{tikzpicture}[overlay, remember picture]
      \node(x)[fit=(a) (b) (b1),inner sep=2pt]{};
      \node(z)[fit=(c),inner sep=2pt]{};
      \filldraw[rounded corners,opacity=.1,red](x.north west)--(x.south west)--(z.south east)--(z.north east)--cycle;
    \end{tikzpicture}
    \caption{Weather-related information set associated with $\mathcal{H}_n$.}
    \label{fig:1}
\end{figure}
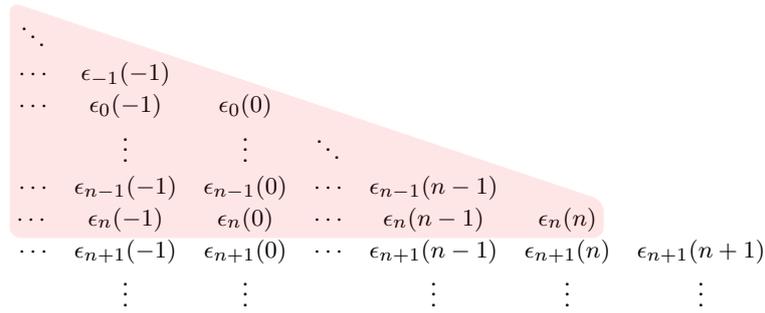

\begin{figure}
    \centering
    \begin{equation*}
    \begin{matrix}
    \tm[a]{\ddots} &                   &                   &        &                     \\
    \cdots & \epsilon_{-1}(-1) &                   &        &                     \\
    \cdots & \epsilon_{0}(-1)  & \epsilon_0(0)     &        &                    \\
           & \vdots            & \vdots            & \ddots &                    \\
    \cdots & \epsilon_{n-1}(-1)  & \epsilon_{n-1}(0) & \cdots & \epsilon_{n-1}(n-1)\\
    \cdots & \epsilon_{n}(-1)  & \epsilon_{n}(0) & \cdots & \epsilon_{n}(n-1)& \tm[d]{\epsilon_{n}(n)}\\
    \cdots & \epsilon_{n+1}(-1)  & \epsilon_{n+1}(0) & \cdots & \epsilon_{n+1}(n-1)& \epsilon_{n+1}(n)\tm[e]{} & \epsilon_{n+1}(n+1)\\
    \tm[b]{} & \vdots  & \vdots & & \vdots& \tm[c]{\vdots} & \vdots
    \end{matrix}    
    \end{equation*}
    
    \begin{tikzpicture}[overlay, remember picture]
      \node(x)[fit=(a) (b),inner sep=2pt]{};
      \node(z)[fit=(c) (e) (d),inner sep=2pt]{};
      \filldraw[rounded corners,opacity=.1,red](x.north west)--(x.south west)--(z.south east)--(z.north east)--cycle;
    \end{tikzpicture}
    \caption{Weather-related information set associated with $\mathcal{G}_n$.}
    \label{fig:2}
\end{figure}
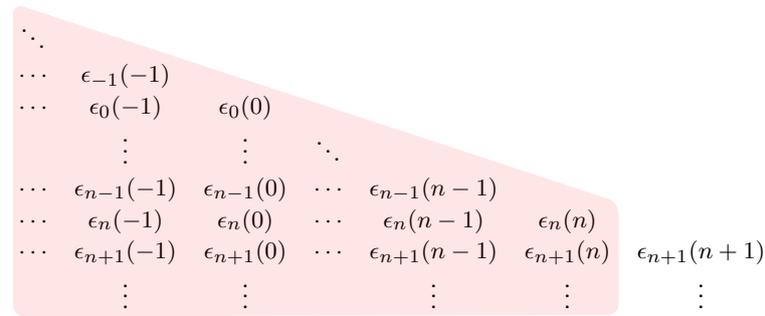

\par We then note that \ref{eq:4.5} implies that for $n\geq k$,
\begin{equation}
\begin{split}
    F_{n|k} &= \mathbb{E}[W_n|\mathcal{G}_k]\\
    &=\mathbb{E}W_0 + G^{n-k}\tilde{W}_k + \sum_{r=-\infty}^{k} \sum_{i=0}^{n-k-1} G^i \epsilon_{n-i}(r),
\end{split}
\label{eq:4.6}
\end{equation}
and the corresponding martingale differences are given by
\begin{equation}
    D_{n|k} = \sum_{i=0}^{n-k} G^i \epsilon_{n-i}(k).
    \label{eq:4.7}
\end{equation}
In addition, we note that \ref{eq:4.6} implies that 
\begin{equation}
    F_{n+1|k}-\mathbb{E}W_0=G\left(G^{n-k}\tilde{W}_k + \sum_{r=-\infty}^{k} \sum_{i=0}^{n-k-1} G^i \epsilon_{n-i}(r)\right) + \sum_{r=-\infty}^{k}\epsilon_{n+1}(r)
\end{equation}
so that
\begin{equation}
    F_{n+1|k}-\mathbb{E}W_0=G(F_{n|k} - \mathbb{E}W_0) + \sum_{r=-\infty}^{k}\epsilon_{n+1}(r).
\end{equation}
Since $(I-G)\mathbb{E}W_0=\mathbb{E}Z_0$, it follows that
\begin{equation}
    F_{n+1|k}=GF_{n|k} + Y_{n+1}(\mathcal{G}_k)
    \label{eq:4.8}
\end{equation}
where
\begin{equation}
    Y_{n+1}(\mathcal{G}_k) \triangleq \mathbb{E}Z_0 + \sum_{r=-\infty}^{k}\epsilon_{n+1}(r).
\end{equation}
As a consequence of \ref{eq:4.8}, we see that the ``forward forecasts'' from time $k$ are correlated, and form (for each $k$) their own state space model with independent (but not identically distributed) ``noise'' rv's $(Y_n(\mathcal{G}_k): n>k)$, initialized at $F_{k|k}=W_k$. Such correlation in the forward forecasts is clearly desirable from a modeling perspective.
\par We now turn to the conditional dynamics of $(W_n: n\geq k)$, conditional on $\mathcal{G}_k$. Define the $W_n(\mathcal{G}_k)$'s via
\begin{equation}
    \mathbb{P}\left((W_n(\mathcal{G}_k): n\geq k)\in\cdot\right) = \mathbb{P}\left((W_n: n\geq k)\in\cdot|\mathcal{G}_k\right).
\end{equation}
The relations \ref{eq:4.5} and \ref{eq:4.6} imply that
\begin{equation}
    \begin{split}
        W_{n+1}(\mathcal{G}_k) - F_{n+1|k} &= \sum_{r=k+1}^{n+1} \sum_{i=0}^{n+1-r} G^i \epsilon_{n+1-i}(r)\\
        &= G\sum_{r=k+1}^{n} \sum_{j=0}^{n-r} G^j \epsilon_{n-j}(r) + \sum_{r=k+1}^{n+1} \epsilon_{n+1}(r)\\
        &= G(W_n(\mathcal{G}_k) - F_{n|k}) + \sum_{r=k+1}^{n+1} \epsilon_{n+1}(r)
    \end{split}
\end{equation}
It follows that for $n\geq k$,
\begin{equation}
    W_{n+1}(\mathcal{G}_k) = GW_{n}(\mathcal{G}_k) + Z_{n+1}(\mathcal{G}_k),
    \label{eq:4.9}
\end{equation}
where
\begin{equation}
    Z_{n+1}(\mathcal{G}_k) = F_{n+1|k} - GF_{n|k}+ \sum_{r=k+1}^{n+1} \epsilon_{n+1}(r).
    \label{eq:27}
\end{equation}

\par Consequently, $(W_n(\mathcal{G}_k): n\geq k)$ is (conditional on $\mathcal{G}_k$) a Markov chain that is a linear state space model driven by a sequence $(Z_n(\mathcal{G}_k): n> k)$ of conditionally independent (but non-identically distributed) rv's. For a given $k$, the variance of the $Z_n(\mathcal{G}_k)$ sequence (conditional on $\mathcal{G}_k$) increases with $n$, so the ``uncertainty plume'' correspondingly grows with time, as one would expect.
\par
\begin{figure}
    \centering
    \begin{equation*}
    \begin{matrix}
    \tm[a]{\ddots} &                   &                   &        &                     \\
    \cdots & \epsilon_{-1}(-1) &                   &        &                     \\
    \cdots & \epsilon_{0}(-1)  & \epsilon_0(0)     &        &                    \\
           & \vdots            & \vdots            & \ddots &                    \\
    \cdots & \epsilon_{n-1}(-1)  & \epsilon_{n-1}(0) & \cdots & \epsilon_{n-1}(n-1)\\
    \cdots & \epsilon_{n}(-1)  & \epsilon_{n}(0) & \cdots & \epsilon_{n}(n-1)& \tm[e]{\epsilon_{n}(n)}\\
    \cdots & \epsilon_{n+1}(-1)  & \tm[d]{\epsilon_{n+1}(0)} & \cdots & \epsilon_{n+1}(n-1)& \epsilon_{n+1}(n) & \epsilon_{n+1}(n+1)\\
    \tm[b]{} & \vdots  & \tm[c]{\vdots} & & \vdots& \vdots & \vdots
    \end{matrix}    
    \end{equation*}
    
    \begin{tikzpicture}[overlay, remember picture]
      \node(x)[fit=(a) (b),inner sep=2pt]{};
      \node(y)[fit=(c) (d),inner sep=2pt]{};
      \node(z)[fit=(e),inner sep=2pt]{};
      \filldraw[rounded corners,opacity=.1,red](x.north west)--(x.south west)--(y.south east)--(y.north east)--(z.south east)--(z.north east)--cycle;
    \end{tikzpicture}
    \caption{Weather-related information set associated with $\mathcal{H}_n\vee\mathcal{G}_0$.}
    \label{fig:3}
\end{figure}
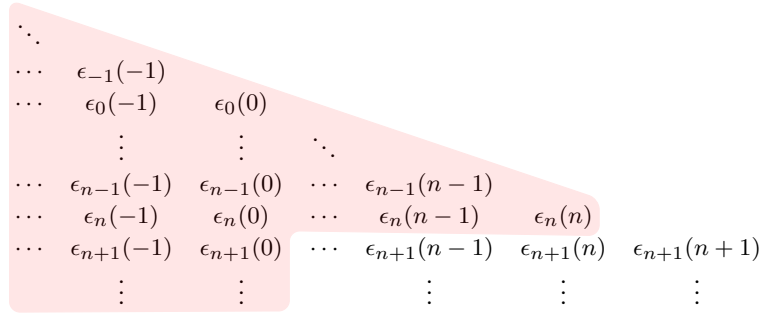

 With \ref{eq:4.9} in hand, we can now modify the value function recursion of Section~\ref{sec:2} so as to compute the optimal policy when the energy system decision maker has available at time $n\geq 0$ the information present in $\mathcal{F}_n\vee\mathcal{K}_0$, the smallest $\sigma$-algebra containing both $\mathcal{F}_n$ and the forecasts collected up to time $0$ by the manager. The structure of our model implies that the policy that is optimal over $\mathcal{F}_n\vee\mathcal{H}_n\vee\mathcal{G}_0$-adapted policies is actually $\mathcal{F}_n\vee\mathcal{K}_0$-measurable. Figure~\ref{fig:3} illustrates the weather-related information set associated with $\mathcal{H}_n\vee\mathcal{G}_0$.
\par Our goal is to minimize
\begin{equation}
    \sum_{j=0}^t \mathbb{E}[c(x_j, A_j, W_{j+1})|\mathcal{G}_0, X_0]
\end{equation}
over all $\mathcal{F}_n\vee\mathcal{H}_n\vee\mathcal{G}_0$-adapted policies $(A_j: 0\leq j \leq t)$. For $1\leq i \leq t$, define the operator $P_{a,i}(\mathcal{G}_0)$ via
\begin{equation}
    (\mathbb{P}_{a,i}(\mathcal{G}_0)h)(x,w)= \int_{\mathbb{R}^{d}}\int_{\mathbb{R}^{m_2}} h(\phi(x, a, Gw + z,v), Gw +z)\mathbb{P}(Z_i(\mathcal{G}_0)\in dz)\mathbb{P}(v_1 \in dv)
\end{equation}
and set
\begin{equation}
    \tilde{c}_i(\mathcal{G}_0, x, a, w) \int_{\mathbb{R}^{d}} c(x, a, Gw+z)\mathbb{P}(Z_i(\mathcal{G}_0)\in dz)
\end{equation}
for $1\leq i \leq t+1$. The appropriate value function backwards recursion in this setting is then given by
\begin{equation}
    v_i(x,w) = \min_{a\in\mathcal{A}}[\tilde{c}_{i+1}(\mathcal{G}_0, x, a, w) + (P_{a,i+1}(\mathcal{G}_0)v_{i+1})(x, w)]
    \label{eq:4.10}
\end{equation}
for $0\leq i< t$, subject to the terminal condition
\begin{equation}
    v_t(x,w) = \min_{a\in\mathcal{A}} \tilde{c}_{t+1}(\mathcal{G}_0, x, a, w).
    \label{eq:4.11}
\end{equation}
As in Section~\ref{sec:2}, an optimal $\mathcal{F}_n\vee\mathcal{H}_n\vee\mathcal{G}_0$-adapted policy is then given by $A_i^* = a_i^*(X_i, W_i)$ for $0\leq i\leq t$, where $a_i^*(x, w)$ is any minimizer of the right-hand side of \ref{eq:4.10} and \ref{eq:4.11}.

\section{MDP's incorporating a dynamic forecast}
\label{sec:5}
In this section, we discuss how our energy system manager should modify her decision-making when she has access to a new set of forecasts each day. More precisely, suppose that at each time $k$ through the decision horizon, the decision maker receives the forecasts $(F_{n|k}: n\geq k)$ prior to making the decision for that period. Now, the decision made at time $k$ can depend on both $W_k$ and the $F_{n|k}$'s. In particular, the decision can now be $\mathcal{F}_k\vee\mathcal{K}_k$-adapted. Since there is more information about $(W_n: n\geq k)$ available when one uses the forecasts, this will typically modify the optimal control relative to the previously discussed formulations of Sections~\ref{sec:2} and ~\ref{sec:4}. Since the forecasts used by the decision maker are constantly updated as $k$ increases, we refer to this setting as a \textit{dynamic forecast} formulation.
\par Let $\Vec{F}_n=(F_{n+j|n}: j\in\mathbb{Z}_+)$ be the entire set of forward forecasts issued at time $n$ (and computed from the history $\mathcal{G}_n$). Recall that $W_n=F_{n|n}$. We claim that the infinite-dimensional process $(\Vec{F}_n : n\in\mathbb{Z})$ is a Markov chain. To see this, observe that \ref{eq:4.7} implies that
\begin{equation}
    \begin{split}
        F_{n+1+j|n+1} &= F_{n+j+1|n} + D_{n+j+1|n+1}\\
         &= F_{n+j+1|n} + \sum_{i=0}^jG^i \epsilon_{n+1+j-i}(n+1)
    \end{split}
    \label{eq:33}
\end{equation}
for $j\geq 0$. Since the collection of rv's $(\epsilon_{n+1+j}(n+1): j\geq 0)$ is independent of $\mathcal{G}_n$, it follows that $(\Vec{F}_n:n\in\mathbb{Z})$ is a Markov chain. One important and related characteristic of our model is that the Markov chain can be initialized from an arbitrary set of values. This means that our model is consistent with any set of forecast values specified at time 
$0$.
\par Of course, we cannot effectively compute optimal policies with a Markov chain having an infinite dimensional state space. So, we need to truncate the set of forecasts that we use within our formulation in order to generate a finite dimensional Markov state variable. In particular, suppose that $\mathcal{G}_{n,r}$ is the smallest $\sigma$-algebra containing both $\mathcal{H}_n$ and the $\sigma$-algebra $\sigma(\epsilon_{n+j}(k):k\leq n, 1\leq j\leq r)$, so that it contains only the forecaster's information about the $r$ future forecasts for periods $n+1, \cdots, n+r$, in addition to the information associated with $\mathcal{H}_n$. Figure~\ref{fig:4} illustrates the weather-related information set associated with $\mathcal{G}_{n,r}$. We note that $\mathcal{G}_n \supseteq \mathcal{G}_{n,r}$ and that for $1\leq j \leq r$, $F_{n+j|n} = \mathbb{E}[W_{n+j}|\mathcal{G}_n]$ is a function only of rv's associated with $\mathcal{G}_{n,r}$, and hence is $\mathcal{G}_{n,r}$-measurable.
\par Using the information associated with $\mathcal{G}_{n,r}$, we can use the recursion in \ref{eq:33} for $0\leq j < r$. For $j=r$, we can use \ref{eq:4.8} to expand $F_{n+1+r|n}$, which yields the recursion
\begin{equation}
\begin{split}
    F_{n+1+r|n+1} &= F_{n+1+r|n} + D_{n+1+r|n+1}\\
    &= GF_{n+r|n} + \mathbb{E}Z_0 + \sum_{j=-\infty}^n \epsilon_{n+1+r}(j) + \sum_{i=0}^r G^i \epsilon_{n+1+r-i}(n+1).
\end{split} 
\label{eq:5.2}
\end{equation}
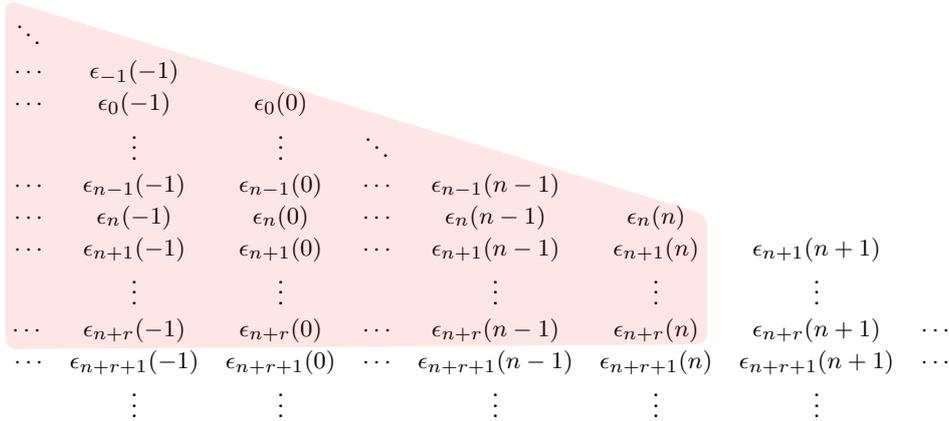
\begin{figure}
    \centering
    \begin{equation*}
    \begin{matrix}
    \tm[a]{\ddots} &                   &                   &        &                     \\
    \cdots & \epsilon_{-1}(-1) &                   &        &                     \\
    \cdots & \epsilon_{0}(-1)  & \epsilon_0(0)     &        &                    \\
           & \vdots            & \vdots            & \ddots &                    \\
    \cdots & \epsilon_{n-1}(-1)  & \epsilon_{n-1}(0) & \cdots & \epsilon_{n-1}(n-1)\\
    \cdots & \epsilon_{n}(-1)  & \epsilon_{n}(0) & \cdots & \epsilon_{n}(n-1)& \epsilon_{n}(n)\tm[d]{}\\
    \cdots & \epsilon_{n+1}(-1)  & \epsilon_{n+1}(0) & \cdots & \epsilon_{n+1}(n-1)& \epsilon_{n+1}(n) & \epsilon_{n+1}(n+1)\\
     & \vdots  & \vdots & & \vdots& \vdots & \vdots\\
    \tm[b]{\cdots} & \tm[b1]{\epsilon_{n+r}(-1)}  & \epsilon_{n+r}(0) & \cdots & \epsilon_{n+r}(n-1)& \epsilon_{n+r}(n)\tm[c]{} & \epsilon_{n+r}(n+1) & \cdots \\
    \cdots & {\epsilon_{n+r+1}(-1)}  & \epsilon_{n+r+1}(0) & \cdots & \epsilon_{n+r+1}(n-1)& \epsilon_{n+r+1}(n) & \epsilon_{n+r+1}(n+1) & \cdots\\
    & \vdots  & \vdots & & \vdots& \vdots & \vdots
    \end{matrix}    
    \end{equation*}
    
    \begin{tikzpicture}[overlay, remember picture]
      \node(x)[fit=(a) (b) (b1),inner sep=2pt]{};
      \node(z)[fit=(c) (d),inner sep=2pt]{};
      \filldraw[rounded corners,opacity=.1,red](x.north west)--(x.south west)--(z.south east)--(z.north east)--cycle;
    \end{tikzpicture}
    \caption{Weather-related information set associated with $\mathcal{G}_{n,r}$.}
    \label{fig:4}
\end{figure}
As a result, $\Vec{F}_{n+1,r}\triangleq(F_{n+1+j|n+1} : 1\leq j\leq r)$ is a linear function of $\Vec{F}_{n,r}$ and a collection of rv's $(\epsilon_{n+1+i}(n+1), \epsilon_{n+1+r}(j) : 1\leq i \leq r, j\leq n)$ that are independent of $\mathcal{G}_{n,r}$. It follows that $(\Vec{F}_{n,r}: n\in\mathbb{Z})$ is an $rd$-dimensional Markov chain. It is also easily seen that it is a Markov chain with stationary transition probabilities. Furthermore, $W_{n+1}=F_{n+1|n+1}$ is a simple stochastic function of $\Vec{F}_{n,r}$, specifically
\begin{equation}
    W_{n+1} = F_{n+1|n} + \epsilon_{n+1}(n+1),
\end{equation}
so that it can easily be generated from $\Vec{F}_{n,r}$ simultaneously with $\Vec{F}_{n+1,r}$.
We can now turn to the computation of the optimal policy in this setting. In particular, we seek the $\mathcal{F}_n\vee\mathcal{G}_{n,r}$-adapted policy that minimizes
\begin{equation}
    \sum_{j=0}^t\mathbb{E}[c(X_j,A_j,W_{j+1})|X_0, \Vec{F}_{0, r}]
\end{equation}
over all $\mathcal{F}_n\vee\mathcal{G}_{n,r}$-adapted policies $(A_n: 0\leq n \leq t)$. Define the operator $P_a$ (acting on integrable functions $h$) via
\begin{equation}
    (P_ah)(x,\Vec{f}) = \mathbb{E}[h(\phi(x, a, W_1, V_1),\Vec{F}_{1,r})|\Vec{F}_{0, r}=\Vec{f})]
\end{equation}
and let
\begin{equation}
    \tilde{c}(x, a, \Vec{f}) = \mathbb{E}[c(x, a, W_1) | \Vec{F}_{0, r}=\Vec{f})].
\end{equation}
We can then compute the associated value functions for this formulation via the backwards recursion
\begin{equation}
    v_i(x,\Vec{f}) = \min_{a\in\mathcal{A}}[\tilde{c}(x, a, \Vec{f}) + (P_{a}v_{i+1})(x, \Vec{f})]
    \label{eq:5.4}
\end{equation}
for $0\leq i< t$, subject to the terminal condition
\begin{equation}
    v_t(x,\Vec{f}) = \min_{a\in\mathcal{A}} \tilde{c}(x, a, \Vec{f}).
    \label{eq:5.5}
\end{equation}
The optimal $\mathcal{F}_n\vee\mathcal{G}_{n,r}$-adapted action $A_n^*$ to be taken in period $n$ is then given by $A_n^*=a_n^*(X_n, \Vec{F}_{n,r})$, where $a_n^*(x,\Vec{f})$ is the minimizer of the right-hand side of \ref{eq:5.4} or \ref{eq:5.5} corresponding to $v_n(x,\Vec{f})$ for $0\leq n \leq t$.
\paragraph{Remark} We note that the use of the reduced Markov state variable $\Vec{F}_{n,r}$  for the weather variables (as opposed to using the state variable $(W_n, \Vec{F}_{n,r})$) is possible only because we made the modeling decision in Section~\ref{sec:2} to express the control state recursion in the form
\begin{equation}
    X_{n+1} = \phi(X_n, A_n, W_{n+1}, V_{n+1})
\end{equation}
and cost $c(X_m, A_n, W_{n+1})$ in terms of $W_{n+1}$ rather than $W_n$. If we had instead modeled the control state evolution via 
\begin{equation}
    X_{n+1} = \phi(X_n, A_n, W_n, V_{n+1})
\end{equation}
and/or cost $c(X_m, A_n, W_n)$, then the decision maker at time $n$ would need to know $W_n$, and $W_n$ would then need to be added to the Markov state variable for the weather. Since either choice, $W_n$ or $W_{n+1}$, is typically reasonable from a modeling viewpoint, we choose to use $W_{n+1}$ in order to obtain this state reduction.

\section{MDP's incorporating both static and dynamic forecasts}
\label{sec:6}
\begin{figure}
    \centering
    \begin{equation*}
    \begin{matrix}
    \tm[a]{\ddots} &                   &                   &        &                     \\
    \cdots & \epsilon_{-1}(-1) &                   &        &                     \\
    \cdots & \epsilon_{0}(-1)  & \epsilon_0(0)     &        &                    \\
           & \vdots            & \vdots            & \ddots &                    \\
    \cdots & \epsilon_{n-1}(-1)  & \epsilon_{n-1}(0) & \cdots & \epsilon_{n-1}(n-1)\\
    \cdots & \epsilon_{n}(-1)  & \epsilon_{n}(0) & \cdots & \epsilon_{n}(n-1)& \epsilon_{n}(n)\tm[f]{}\\
    \cdots & \epsilon_{n+1}(-1)  & \epsilon_{n+1}(0) & \cdots & \epsilon_{n+1}(n-1)& \epsilon_{n+1}(n) & \epsilon_{n+1}(n+1)\\
     & \vdots  & \vdots & & \vdots& \vdots & \vdots\\
    \cdots & \epsilon_{n+r}(-1)  & \epsilon_{n+r}(0) & \cdots & \epsilon_{n+r}(n-1)& \epsilon_{n+r}(n)\tm[e]{} & \epsilon_{n+r}(n+1) & \cdots \\
    \cdots & \epsilon_{n+r+1}(-1)  & \tm[d]{\epsilon_{n+r+1}(0)} & \cdots & \epsilon_{n+r+1}(n-1)& \epsilon_{n+r+1}(n)\tm[c]{} & \epsilon_{n+r+1}(n+1) & \cdots \\
    \tm[b]{}& \vdots  & \vdots\tm[c]{} & & \vdots& \vdots & \vdots
    \end{matrix}    
    \end{equation*}
    
    \begin{tikzpicture}[overlay, remember picture]
      \node(x)[fit=(a) (b),inner sep=2pt]{};
      \node(y)[fit=(c) (d),inner sep=2pt]{};
      \node(z)[fit=(e) (f),inner sep=2pt]{};
      \filldraw[rounded corners,opacity=.1,red](x.north west)--(x.south west)--(y.south east)--(y.north east)--(z.south east)--(z.north east)--cycle;
    \end{tikzpicture}
    \caption{Weather-related information set associated with $\mathcal{G}_0\vee\mathcal{G}_{n, r}$.}
    \label{fig:5}
\end{figure}
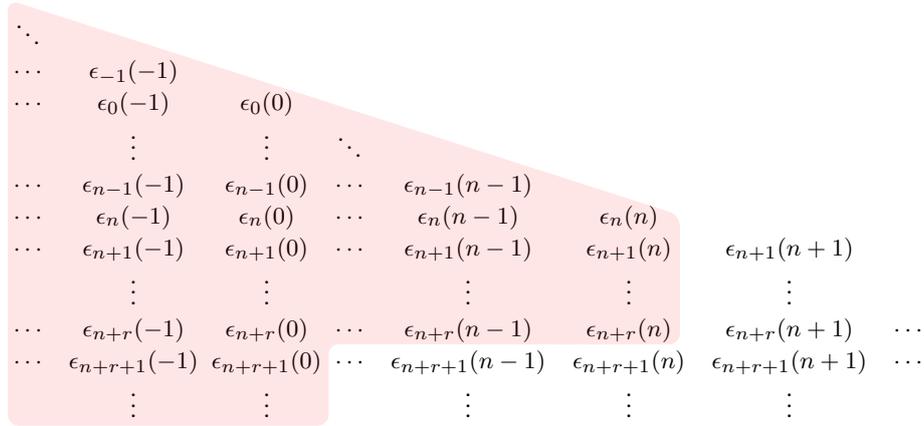

For computational tractability, the value of $r$ used in $\Vec{F}_{n,r}$ will typically need to be small. But weather forecasters will typically provide forward forecasts over a much larger number of periods. In order to (partially) account for these longer range forecasts (without expanding our state description for the MDP), we now build a formulation that takes into account all the forward forecasts that are present in $\mathcal{G}_0$ (i.e. the static forecasts that are available at time $0$), as well as the dynamic forecasts associated with $\mathcal{G}_{n,r}$ for $1\leq n \leq t$. Thus, in this formulation, the decision maker at time $n$ has access to $X_n, F_{n+1|n}, \cdots,F_{n+r|n}$ and $F_{j|0}$ for $j\geq 1$. Figure~\ref{fig:5} illustrates the weather-related information set corresponding to $\mathcal{G}_0\vee\mathcal{G}_{n,r}$.
\par We now turn to the conditional dynamics of $(\Vec{F}_{n,r}:n\geq 0)$, conditional on $\mathcal{G}_0$. Define the $\Vec{F}_{n,r}(\mathcal{G}_0)$'s via
\begin{equation}
    \mathbb{P}((\Vec{F}_{n,r}(\mathcal{G}_0): n\geq 0)\in\cdot) = \mathbb{P}((\Vec{F}_{n,r}: n\geq 0)\in\cdot|\mathcal{G}_0)
\end{equation}
Because the martingale difference terms on the right-hand side of \ref{eq:33} are independent of $\mathcal{G}_0$,
\begin{equation}
 F_{n+1+j|n+1}(\mathcal{G}_0) = F_{n+1+j|n}(\mathcal{G}_0) + \beta_{n+1,j}(\mathcal{G}_0)
 \label{eq:6.2}
\end{equation}
for $0\leq j < r$, where $\beta_{n+1,j}(\mathcal{G}_0)$ is independent of $\mathcal{G}_0\vee\mathcal{G}_{n,r}$ and
\begin{equation}
 \beta_{n+1,j}(\mathcal{G}_0) \stackrel{D}{=} \sum_{i=0}^jG^i\epsilon_{n+1+j-i}(n+1)
\end{equation}
On the other hand, the right-hand side of \ref{eq:5.2} contains terms that are $\mathcal{G}_0$-measurable. In particular,
\begin{equation}
 F_{n+1+r|n+1}(\mathcal{G}_0) = GF_{n+r|n}(\mathcal{G}_0) + \mathbb{E}Z_0+ \beta_{n+1,r}(\mathcal{G}_0)
 \label{eq:47}
\end{equation}
where
\begin{equation}
 \beta_{n+1,r}(\mathcal{G}_0) = \sum_{j=-\infty}^0 \epsilon_{n+1+r}(j) + \Lambda_{n+1, r}
 \label{eq:6.5}
\end{equation}
where $\sum_{j=-\infty}^0 \epsilon_{n+1+r}(j)$ is $\mathcal{G}_0$-measurable and $\Lambda_{n+1, r}$ is independent of $\mathcal{G}_0\vee\mathcal{G}_{n,r}$, where
\begin{equation}
 \Lambda_{n+1, r}\stackrel{D}{=} \sum_{j=1}^n \epsilon_{n+1+r}(j) + \sum_{i=0}^r G^i \epsilon_{n+1+r-i}(n+1).
\end{equation}
Since $\Vec{F}_{n+1, r}(\mathcal{G}_0)$ can be expressed as a function of $\Vec{F}_{n, r}(\mathcal{G}_0)$ and a family of rv's $\beta_{n+1}(\mathcal{G}_0) \triangleq (\beta_{n+1,j}(\mathcal{G}_0) : 1\leq j \leq r)$ that are independent of $\mathcal{H}_n^0$, it follows that $(\Vec{F}_{n, r}(\mathcal{G}_0) : n\geq 0)$ is (conditional on $\mathcal{G}_0$) a Markov chain. However, as with the Markov chain of Section~\ref{sec:4}, the conditioning on $\mathcal{G}_0$ makes this a Markov chain with non-stationary transition probabilities; see \ref{eq:6.5} in particular. Furthermore, as in Section~\ref{sec:4}, the variance of $\beta_{n+1,r}(\mathcal{G}_0)$ (conditional on $\mathcal{G}_0$) increases in $n$, so that the ``uncertainty plume'' increases over time.
\par We now turn to the computation of a policy $(A_n^*: 0\leq n \leq t)$ that minimizes
\begin{equation}
    \mathbb{E}[\sum_{j=0}^tc(X_j,A_j,W_{j+1})|X_0,\mathcal{G}_0]
\end{equation}
over all policies $(A_n: 0\leq n \leq t)$ that are $\mathcal{H}_n^0$-adapted. For $\Vec{f}=(f_1, f_2, \cdots, f_r)$, define the operator $P_{a,i}(\mathcal{G}_0)$ (acting on integrable functions $h$) via
\begin{equation}
    \begin{split}
    (P_{a,i}(\mathcal{G}_0) &h)(x, \Vec{f})= \int_{\mathbb{R}^{(r+1) d}} \int_{\mathbb{R}^{m_2}} h(\phi(x, a, w, v), y_1, y_2, \dots, y_r) \mathbb{P}(f_1 + \beta_{i, 0}(\mathcal{G}_0)\in dw, \\
    &f_{j+1} + \beta_{i, j}(\mathcal{G}_0)\in dy_j, 1\leq j\leq r, Gf_r + \mathbb{E}Z_0 + \beta_{i, r}(\mathcal{G}_0)\in dy_r)\mathbb{P}(V_1 \in dv)
    \end{split}
\end{equation}
and set
\begin{equation}
    \tilde{c}_i(\mathcal{G}_0, x, a, \Vec{f}) \int_{\mathbb{R}^{d}} c(x, a, w)\mathbb{P}(f_1 + \beta_{i,0}(\mathcal{G}_0)\in dw)
\end{equation}
for $1\leq i\leq t+1$. As in Section~\ref{sec:5}, the value function recursion takes the form
\begin{equation}
    v_i(x,\Vec{f})=\min_a[\tilde{c}_{i+1}(\mathcal{G}_0, x, a, \Vec{f}) + (P_{a,i}(\mathcal{G}_0)v_{i+1})(x, \Vec{f})]
    \label{eq:6.10}
\end{equation}
 for $1\leq i< t$, with terminal condition
\begin{equation}
    v_{t}(x, \Vec{f}) = \min_a \tilde{c}_t(\mathcal{G}_0, x, a, \Vec{f}).
    \label{eq:6.11}
 \end{equation}
Again, the optimal $\mathcal{H}_{n}^0$-adapted policy is then given by $A_i^*=a_i^*(X_i, \Vec{F}_i)$ for $0\leq i \leq t$, where $a_i^*(x, \Vec{f})$ is any minimizer of the right-hand side of \ref{eq:6.10} and \ref{eq:6.11}.
 
\section{An energy control system example}
\label{sec:7}
In this section, we illustrate some of our theory in the setting of a simple energy control system example. In particular, we let $W_n$ represent the ambient outdoors temperature at the beginning of period $n$ at the site of the energy system that is under control. We assume that $(W_n:n\in\mathbb{Z})$ is a real-valued Markov chain corresponding to a first order autoregressive process, so that
\begin{equation}
    W_{n+1}=gW_n+Z_{n+1}
    \label{eq:7.1}
\end{equation}
for $n\in\mathbb{Z}$, where $g\in(0,1)$ and the $Z_i$'s are iid with $\mathbb{E}Z_0^2<\infty$. To help interpret $g$, we note that $\text{corr}(W_j, W_{j+n})=g^n$, so that the number of periods for the correlation to decay to $0.1$ is approximately $\log(0.1)/\log g$.
\par We now describe our simplified energy control system corresponding to heating and cooling a building. We assume that the difference $\Delta_n(\triangleq X_n-W_n)$ between the internal ($X_n$) and external ($W_n$) temperatures is ``mean reverting'', so that the $\Delta_n$'s satisfy their own first order autoregression. In particular, in the absence of control,
\begin{equation}
    \Delta_{n+1}=\rho\Delta_n+V_{n+1},
    \label{eq:7.2}
\end{equation}
where the $V_j$'s are iid and independent of the $Z_k$'s with $\mathbb{E}V_0^2<\infty$. We expect the building to equilibrate more rapidly than does the outdoors temperature, so we expect $\rho\in(0, g)$. Substituting \ref{eq:7.1} into \ref{eq:7.2}, we find that in the presence of the control $A_n$,
\begin{equation}
    X_{n+1}=(g-\rho)W_n + \rho X_n + Z_{n+1} + V_{n+1} + A_n
    \label{eq:7.3}
\end{equation}
for $n\in\mathbb{Z}$.
\par We now wish to take advantage of the powerful toolset that is available when our state space model has a quadratic cost structure. We assume that our goal is to minimize the expected infinite horizon discounted cost given by
\begin{equation}
    \mathbb{E}\sum_{j=0}^{\infty}\alpha^j\left[(X_j-\tau)^2+\kappa A_j^2\right]
    \label{eq:7.4}
\end{equation}
over all $\mathcal{F}_j$-adapted controls, where $\kappa>0$, $\alpha \in (0,1)$ is the discount factor, and $\tau$ is the reference temperature to which we are trying to steer the system. To incorporate $\tau$ into the linear/quadratic formulation, we add $Y_j$ as a state variable for which $Y_j=Y_{j-1}$ for $j\in \mathbb{Z}$. Furthermore, we let $\tilde{W}_j=W_j-(1-g)^{-1}\mathbb{E}Z_0$, $\tilde{Z}_j=Z_j-\mathbb{E}Z_0$, $\tilde{X}_j=X_j-\tau_0$, $\tau_0=(1-g)^{-1}\mathbb{E}Z_0+(1-\rho)^{-1}\mathbb{E}V_0$ and $\tilde{V}_j=V_j-\mathbb{E}V_0$, and rewrite \ref{eq:7.1} and \ref{eq:7.3} in terms of the mean zero ``noise'' rv's $\tilde{Z}_{n+1}$ and $\tilde{V}_{n+1}$:
\begin{equation}
\begin{split}
    \tilde{W}_{n+1} &= g\tilde{W}_n + \tilde{Z}_{n+1}\\
    \tilde{X}_{n+1} &= (g-\rho)\tilde{W}_n + \rho\tilde{X}_n + \tilde{Z}_{n+1} + \tilde{V}_{n+1} + A_n.
\end{split}
\end{equation}
Furthermore, we can express $X_n-\tau$ as $\tilde{X}_n-Y_n$, where we take $Y_0=\tau -\tau_0$.
\par Set $\chi_n=(\tilde{W}_n, \tilde{X}_n, Y_n)^T$ and $\xi_n=(\tilde{Z}_n, \tilde{Z}_n + \tilde{V}_n, 0)^T$, and observe that we can express our control system dynamics as
\begin{equation}
    \chi_{n+1}=A\chi_n+BA_n+\xi_{n+1}
    \label{eq:7.6}
\end{equation}
for $n\geq 0$, where
\begin{equation}
    A=\begin{pmatrix} g & 0 & 0\\g-\rho & \rho & 0\\0 & 0 & 1 \end{pmatrix}, \qquad B=\begin{pmatrix} 0\\1\\0\end{pmatrix}.
\end{equation}
The objective \ref{eq:7.4} can then be re-expressed as
\begin{equation}
    \mathbb{E}\sum_{j=0}^{\infty}\alpha^j[\chi_j^TQ\chi_j  + A_j^TRA_j],
    \label{eq:7.8}
\end{equation}
where 
\begin{equation}
    Q=\begin{pmatrix} 0 & 0 & 0\\0 & 1 & -1\\0 & -1 & 1 \end{pmatrix}, \qquad R=\kappa.
\end{equation}
We observe that this model does not satisfy the standard controllabillity hypothesis that is commonly used within the literature on state space models with quadratic costs (in particular, the $Y_j$'s are not controllable). Nevertheless, the special problem structure here allows us to follow the approach on p. 231-233 of \cite{bertsekas2012dynamic} to obtain the solution of this stochastic control problem in closed form.
\par In particular, define the optimal return operator $\mathcal{T}$ (defined on suitably integrable functions $h$) via
\begin{equation}
    (\mathcal{T}h)(\Vec{z}) = \min_a[\Vec{z}^TQ\Vec{z} + a^2\kappa + \alpha\mathbb{E}h(A\Vec{z}+Ba+\xi_1)]
\end{equation}
for $\Vec{z}=(w, x, y)^T$ and note that this stochastic control problem corresponds to a positive dynamic program; see p. 214 of \cite{bertsekas2012dynamic}. Hence, if $v_0(\Vec{x})\equiv 0$, it follows that $\mathcal{T}v_0\geq v_0$, thereby implying that $v_k\equiv \mathcal{T}^kv_0\geq \mathcal{T}^{k-1}v_0=v_{k-1}$, so that the $v_k$'s are functions that are monotone increasing in $k$. Hence,
\begin{equation}
    v_{\infty} = \lim_{k\to\infty}v_k
\end{equation}
exists. Furthermore, if $h_0(\Vec{z})\equiv 0$ and
\begin{equation}
    h_{k+1}(\Vec{z}) = \Vec{z}^T Q \Vec{z} + \alpha\mathbb{E}h_{k}(A\Vec{z}+\xi)
\end{equation}
for $k\geq 0$, $v_k\leq h_k$ for $k\geq 0$. As for the $v_k$'s, $(h_k : k\geq 0)$ is also a monotone sequence so that $h_k\to h_{\infty}$. The limit $h_{\infty}$ is the value function corresponding to the policy in which $A_k=0$ for $k\geq 0$. Since $|g|<1$ and $|\rho|<1$, the associated stochastic dynamical system is stable and $h_{\infty}$ is finite-valued. We can therefore conclude that $v_{\infty}$ is finite-valued.
\par We further note that if $J=(J(i,k): 1\leq i, k\leq 3)$ is a symmetric non-negative definite matrix, the scalar $\alpha B^TJB +R = \alpha J(2,2) + \kappa>0$ (since the diagonal entries of such a matrix must be non-negative), so that $\alpha B^TJB +R$ is guaranteed to be non-singular. As a result, the matrix recursion
\begin{equation}
    K_{j+1}=A^T(\alpha K_j-\alpha^2K_jB(\alpha B^TK_j B+R)^{-1}B^TK_j)A+Q,
    \label{eq:ricatti}
\end{equation}
subject to $K_0=Q$, is well-defined, and
\begin{equation}
    v_{j+1}(\Vec{z}) = \Vec{z}^TK_j\Vec{z} + \sum_{i=0}^{j-1}\alpha^{j-i}\mathbb{E} \xi_1^TK_i\xi_1;
\end{equation}
see p. 231 of \cite{bertsekas2012dynamic}. By following the argument on p. 156 of \cite{bertsekas2005dynamic}, we can conclude that there exists a finite-valued non-negative definite matrix $K_{\infty}=(K_{\infty}(i,k):1\leq i,k\leq 3)$ for which $K_j\to K_\infty$ as $j\to \infty$. Taking limits in \ref{eq:ricatti}, we find that $K_{\infty}$ satisfies the matrix Ricatti equation
\begin{equation}
    K_{\infty}=A^T(\alpha K_{\infty}-\alpha^2K_{\infty}B(\alpha B^TK_\infty B+R)^{-1}B^TK_\infty)A + Q.
\end{equation}
Furthermore, as seen from p. 232 of \cite{bertsekas2012dynamic}, we conclude that the optimal value function for the control problem is
\begin{equation}
    v_{\infty}(\Vec{z}) = \Vec{z}^TK_\infty \Vec{z} + \frac{\alpha}{1 - \alpha}\mathbb{E}\xi_1^T K_\infty \xi_1,
\end{equation}
and the associated optimal action $A^*_j$ to be taken at time $j$ is
\begin{equation}
    \begin{split}
        A^*_j &= -\alpha(\alpha B^TK_\infty B+\kappa)^{-1} B^T K_{\infty}A\chi_j\\
        &= -\alpha(\alpha K_\infty(2,2) +\kappa)^{-1}[ K_\infty(2,1)g\tilde{W}_j + K_\infty(2,2)((g-\rho)W_j + \rho \tilde{X}_j)\\
        &+ K_\infty(3,3)[\tau - \tau_0]].
    \end{split}
\end{equation}
We now turn to the analyzing exactly the same MDP when the dynamic forecasts of Section~\ref{sec:5} are incorporated into the problem. To simplify our exposition, we  set $r=2$, so that our energy system manager has access to the forecasts $F_{n+1|n}$ and $F_{n+2|n}$ (in addition to $W_n$ and $X_n$) at the time that the decision at time $n$ is taken. Put $\tilde{F}_{n+i|n}=F_{n+i|n}-\mathbb{E}Z_0/(1-g)$, for $i=1,2$ and $n\in\mathbb{Z}$. If $\underaccent{\tilde}{\chi}_n = (\tilde{W}_n, \tilde{F}_{n+1|n}, \tilde{F}_{n+2|n}, \tilde{X}_n, Y_n)^T$, Section~\ref{sec:5}'s discussion establishes that
\begin{equation}
    \underaccent{\tilde}{\chi}_{n+1} = \mathcal{A}\underaccent{\tilde}{\chi}_n + \mathcal{B}A_n + \underaccent{\tilde}{\xi}_{n+1},
    \label{eq:7.9}
\end{equation}
where
\begin{equation}
    \mathcal{A}=\begin{pmatrix}
    0 & 1 & 0 & 0 & 0 \\
    0 & 0 & 1 & 0 & 0\\
    0 & 0 & g & 0 & 0\\
    -\rho & 1 & 0 & \rho & 0\\
    0 & 0 & 0 & 0 & 1 \end{pmatrix}, \qquad \mathcal{B}=\begin{pmatrix} 0\\0\\0\\1\\0\end{pmatrix},
\end{equation}
and
\begin{equation}
    \underaccent{\tilde}{\xi}_{n}=\begin{pmatrix} \epsilon_{n+1}(n+1)\\ \epsilon_{n+2}(n+1) + g \epsilon_{n+1}(n+1) \\ \epsilon_{n+3}(n+1) + g \epsilon_{n+2}(n+1) + g^2 \epsilon_{n+1}(n+1) + \sum_{j=-\infty}^n \epsilon_{n+3}(j) \\  \epsilon_{n+1}(n+1) + \tilde{V}_{n+1}\\0\end{pmatrix}.
\end{equation}
We can also re-express the objective \ref{eq:7.4} in terms of the $\underaccent{\tilde}{\chi}_{j}$'s, namely
\begin{equation}
    \mathbb{E}\sum_{j=0}^{\infty}\alpha^j \left[\underaccent{\tilde}{\chi}_j^T \mathcal{Q} \underaccent{\tilde}{\chi}_j  + A_jRA_j\right]
    \label{eq:7.10}
\end{equation}
where
\begin{equation}
    \mathcal{Q}=\begin{pmatrix} 0 & 0 & 0 & 0 & 0 \\0 & 0 & 0 & 0 & 0\\ 0 & 0 & 0 & 0 & 0\\0 & 0 & 0 & 1 & -1\\ 0 & 0 & 0 & -1 & 1 \end{pmatrix}.
\end{equation}
We now assume a specific probability structure on the $\epsilon_{n+j}(n)$'s and $V_n$'s, namely that the $\epsilon_{n+j}(n)$'s and $V_n$'s are independent normally distributed rv's with $\text{var}\, \epsilon_{n+j}(n)=\sigma^2\gamma^{2j}$ and $\text{var}\, V_n\triangleq \sigma_V^2$, where $\gamma\in(0,1)$. With this choice, the $\underaccent{\tilde}{\xi}_{n}$'s are iid multivariate normally distributed random vectors with common matrix covariance $C$ defined by
\begin{equation}
    C=\begin{pmatrix} \sigma^2 & \sigma^2 g & \sigma^2 g^2 & \sigma^2 & 0 \\ \sigma^2 g & \sigma^2 (\gamma^2 + g^2) & \sigma^2  (\gamma^2g + g^3) & \sigma^2 g & 0\\ \sigma^2 g^2 & \sigma^2  (\gamma^2g + g^3) & \sigma^2  (\gamma^4 + \gamma^2g^2 + g^4 + \frac{\gamma^6}{1-\gamma^2}) & \sigma^2 g^2 & 0\\ \sigma^2 & \sigma^2g & \sigma^2g^2 & \sigma^2+ \sigma_V^2 & 0\\ 0 & 0 & 0 & 0 & 0 \end{pmatrix}.
\end{equation}
The solution of the stochastic control problem defined by the quadratic objective \ref{eq:7.10} subject to the state specific dynamics \ref{eq:7.9} follows the same lines as for the control system \ref{eq:7.6} and \ref{eq:7.8} without forecasts. In particular, there exists a solution to the matrix Ricatti equation
\begin{equation}
    \mathcal{K}_{\infty}=\mathcal{A}^T(\alpha \mathcal{K}_{\infty}-\alpha^2\mathcal{K}_{\infty}\mathcal{B}(\alpha \mathcal{K}_{\infty}(4,4)+\kappa)^{-1}\mathcal{B}^T\mathcal{K}_\infty)\mathcal{A} + \mathcal{Q}
\end{equation}
where $\mathcal{K}_{\infty}(4,4)$ is the $(4,4)$ entry of $\mathcal{K}_{\infty}$ and the optimal value function $\underaccent{\tilde}{v}_{n}$ for this MDP is given by
\begin{equation}
    \underaccent{\tilde}{v}_{\infty}(\underaccent{\tilde}{\Vec{z}}) = \underaccent{\tilde}{\Vec{z}}^T \mathcal{K}_{\infty} \underaccent{\tilde}{\Vec{z}} - \frac{\alpha}{1-\alpha}\mathbb{E} \underaccent{\tilde}{\xi}^T \mathcal{K}_{\infty} \underaccent{\tilde}{\xi}
\end{equation}
where $\underaccent{\tilde}{\Vec{z}}=(w,f_1, f_2, x, y)^T$. Furthermore, the optimal action adapted to $\mathcal{G}_{n, r}$ is given by
\begin{equation}
    A^*_n = -\alpha (\alpha \mathcal{K}_{\infty}(4,4)+\kappa)^{-1}\mathcal{B}\mathcal{K}_{\infty}\mathcal{A}\underaccent{\tilde}{\chi}_{n}.
\end{equation}
To compare the two optimal controls, we note that the system is controlled only over $[0, \infty)$ and uncontrolled over $(-\infty, 0)$. Hence, the distribution of the control system at time $0$ is given by the distribution of $((\tilde{W}_0, \tilde{X}_0):n \leq 0)$ for the system without forecasts and by the distribution of $(\underaccent{\tilde}{\chi}_{0}:n\leq 0)$ for the system with forecasts.
\par The expected cost for the system with the optimal $\mathcal{F}_n$-adapted policy is therefore $\mathbb{E}v_\infty(\tilde{W}_0, \tilde{X}_0, \tau - \tau_0)$, while the expected cost for the optimal $\mathcal{G}_{n,r}$-adapted policy is $\mathbb{E}\underaccent{\tilde}{v}_{\infty}(\underaccent{\tilde}{\chi}_{0})$. To explicitly compute the distribution of $(\tilde{W}_0, \tilde{X}_0)$, we note that the uncontrolled system is Gaussian with zero mean, so the distribution of $(\tilde{W}_0, \tilde{X}_0)$ is completely determined by its covariance structure. Because the $W_j$'s and $\Delta_j$'s are indepen\textbf{}dent first order autoregressive processes, $\text{var}\, \tilde{W}_0=(1-g^2)^{-1}\text{var}\, Z_0$ and $\text{var}\, \tilde{\Delta}_0=(1-\rho^2)^{-1}\text{var}\, V_0$. Since $X_0= W_0 + \Delta_0$, it follows that $\text{var}\, \tilde{X}_0 = (1-g^2)^{-1}\text{var}\, Z_0 + (1-\rho^2)^{-1}\text{var}\, V_0$ and $\text{cov}\, (\tilde{W}_0, \tilde{X}_0)= (1-g^2)^{-1}\text{var}\, Z_0$.
\par It follows that
\begin{equation}
\begin{split}
    \mathbb{E}v_\infty&\left(\tilde{W}_0, \tilde{X}_0, \tau - (1-g)^{-1}\mathbb{E}Z_0 - (1 - \rho)^{-1}\mathbb{E}V_0\right)\\
    &= K_{\infty}(1,1)\frac{\text{var}\, Z_0}{1-g^2} + 2K_{\infty}(1,2)\frac{\text{var}\, Z_0}{1-g^2} + K_{\infty}(2,2)\left(\frac{\text{var}\, Z_0}{1-g^2} + \frac{\text{var}\, V_0}{1-\rho^2}\right)\\
    &+ K_{\infty}(3,3) \left(\tau - (1-g)^{-1}\mathbb{E}Z_0 - (1 - \rho)^{-1}\mathbb{E}V_0\right)^2\\
    &+ \frac{\alpha}{1-\alpha}\left[K_{\infty}(1,1)\text{var}\, Z_0 + 2K_{\infty}(1,2)\text{var}\, Z_0 + K_{\infty}(2,2)(\text{var}\, Z_0 + \text{var}\, V_0)\right]
    \label{eq:7.27}
\end{split}
\end{equation}
\par We turn next to the evaluation of $\mathbb{E}\underaccent{\tilde}{v}_{\infty}(\underaccent{\tilde}{\chi}_{0})$. We write $\underaccent{\tilde}{\chi}_{n}^T=(\underaccent{\tilde}{\chi}_{n}(1), \cdots, \underaccent{\tilde}{\chi}_{n}(5))$ and $\underaccent{\tilde}{\xi}_{n}^T=(\underaccent{\tilde}{\xi}_{n}(1), \cdots, \underaccent{\tilde}{\xi}_{n}(5))$. With this notation in hand, we can write
\begin{equation}
    \mathbb{E}\underaccent{\tilde}{v}_{\infty}(\underaccent{\tilde}{\chi}_{0}) = \sum_{i,j=1}^5 \mathcal{K}_\infty(i,j)[\mathbb{E}\underaccent{\tilde}{\chi}_{0}(i)\underaccent{\tilde}{\chi}_{0}(j) + C(i,j)],
    \label{eq:7.28}
\end{equation}
where $C(i,j)$ is the $(i,j)$'th entry of the covariance matrix $C$. As for $\mathbb{E}\underaccent{\tilde}{\chi}_{0}(i)\underaccent{\tilde}{\chi}_{0}(j)$, note that $\mathbb{E}\underaccent{\tilde}{\chi}_{0}(i)\underaccent{\tilde}{\chi}_{0}(5)=0$ for $1\leq i\leq4$, and
\begin{equation}
    \mathbb{E}\underaccent{\tilde}{\chi}_{0}(5)^2=\left(\tau-\frac{\mathbb{E}Z_0}{1-g}-\frac{\mathbb{E}V_0}{1-\rho}\right)^2.
\end{equation}
\par We again note that $X_0=W_0+\Delta_0$, so that $\mathbb{E}\underaccent{\tilde}{\chi}_{0}(1)^2=\text{var}\,\tilde{W}_0 = (1-g^2)^{-1}\text{var}\,Z_0$, $\mathbb{E}\underaccent{\tilde}{\chi}_{0}(4)^2=\text{var}\,\tilde{X}_0 = (1-g^2)^{-1}\text{var}\,Z_0 + (1-\rho^2)^{-1}\text{var}\,V_0$ and $\mathbb{E}\underaccent{\tilde}{\chi}_{0}(i)\underaccent{\tilde}{\chi}_{0}(4)= \mathbb{E}\underaccent{\tilde}{\chi}_{0}(i)(W_0 + \Delta_0)=\mathbb{E}\underaccent{\tilde}{\chi}_{0}(i)W_0 = \mathbb{E}\underaccent{\tilde}{\chi}_{0}(i)\underaccent{\tilde}{\chi}_{0}(1)$ for $1\leq i\leq 3$. Also, since $(\tilde{F}_{j+2|j}:j\leq 0)$ is a first order autoregressive process, $\mathbb{E}\underaccent{\tilde}{\chi}_{0}^2(3) = \text{var}\,\tilde{F}_{2|0}=C(3,3)/(1-g^2)$. To compute a closed form for the remaining entries of $(\mathbb{E}\underaccent{\tilde}{\chi}_{0}(i)\underaccent{\tilde}{\chi}_{0}(j): 1\leq i,j\leq 5)$, we note that
\begin{equation}
    \underaccent{\tilde}{\chi}_{0}\stackrel{\mathcal{D}}{=} \mathcal{A} \underaccent{\tilde}{\chi}_{0} + \underaccent{\tilde}{\xi}_{1}.
\end{equation}
For example, we have $\underaccent{\tilde}{\chi}_{0}(2) \stackrel{\mathcal{D}}{=} \underaccent{\tilde}{\chi}_{0}(3) + \underaccent{\tilde}{\xi}_{0}(2)$, so that
\begin{equation}
    \mathbb{E}\underaccent{\tilde}{\chi}_{0}(2)^2 = \mathbb{E}\underaccent{\tilde}{\chi}_{0}(3)^2 + C(2,2) = \frac{C(3,3)}{1-g^2} + C(2,2).
\end{equation}
Similarly, we find that
\begin{equation}
    \begin{split}
        \mathbb{E}\underaccent{\tilde}{\chi}_{0}(1)\underaccent{\tilde}{\chi}_{0}(2) &= \frac{gC(3,3)}{1-g^2} + C(2,3) + C(1,2),\\
        \mathbb{E}\underaccent{\tilde}{\chi}_{0}(1)\underaccent{\tilde}{\chi}_{0}(3) &= \frac{g^2C(3,3)}{1-g^2} + C(2,3) + C(1,3),\\
        \mathbb{E}\underaccent{\tilde}{\chi}_{0}(2)\underaccent{\tilde}{\chi}_{0}(3) &= \frac{gC(3,3)}{1-g^2} + C(2,3).
    \end{split}
\end{equation}
\par As an alternative to using the closed forms \ref{eq:7.27} and \ref{eq:7.28} to compute $\mathbb{E}v_\infty(\tilde{W}_0, \tilde{X}_0, \tau -\tau_0)$ and $\mathbb{E}\underaccent{\tilde}{v}_{\infty}(\underaccent{\tilde}{\chi}_{0})$, an iterative approach can be used to compute the covariance matrices; this may be preferable for some applications. Given the recursion $\chi_{n+1}=A \chi_{n} + \xi_{n+1}$,
\begin{equation}
    \chi_n \stackrel{\mathcal{D}}{=} \sum_{k=0}^{\infty} A^k \xi_{k}.
\end{equation}
It follows that the covariance matrix for $\chi_{n}$ satisfies
\begin{equation}
\Lambda=A\Lambda A^T+\Sigma_\xi,
\end{equation}
subject to $\Lambda_0=0$ \cite{dejong2021}.
\par
\begin{figure}[tbhp]
\centering\includegraphics[width=.9\columnwidth]{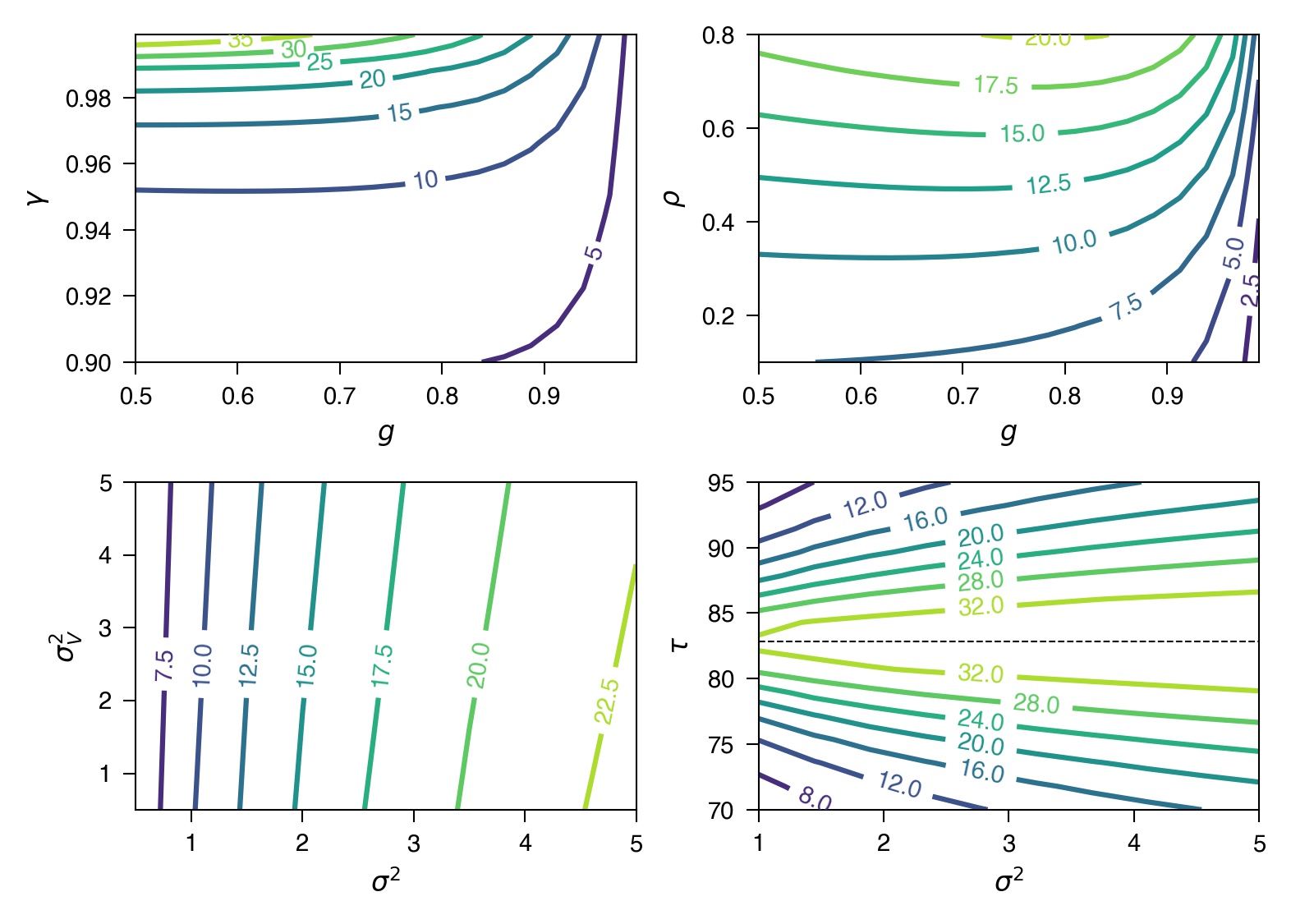}
\caption{Relative improvement $D$ (\%) in \ref{eq:7.33} from using dynamic forecasts versus no forecasts for the energy system control example of Section~\ref{sec:7}. Contour plots are shown for $D$ as a function of pairs of system parameters. In the bottom right plot, a dashed line is drawn at $\tau_0=(1-g)^{-1}\mathbb{E}Z_0 + (1 - \rho)^{-1}\mathbb{E}V_0$. Unless otherwise specified, values used for parameters are $\mathbb{E}W_0=80 F$, $\mathbb{E}V_0=2 F$, $\mathbb{E}W_0=80 F$, $\gamma=.95$, $\alpha=.9$, $g=.6$, $\rho=.3$, $\tau=74 F$, $\sigma^2=\sigma_V^2=1$.}
\label{fig:7.1}
\end{figure}

 Implementations for both the closed-form and the iterative approaches are provided in the case of our energy control system example in the online repository for this paper \cite{dechalendar2021}. Figure~\ref{fig:7.1} shows the percentage reduction in cost from using dynamic forecasts, namely
\begin{equation}
    D \triangleq  \frac{\mathbb{E}v_\infty\left(\tilde{W}_0, \tilde{X}_0, \tau - \tau_0\right) - \mathbb{E}\underaccent{\tilde}{v}_{\infty}(\underaccent{\tilde}{\chi}_{0})}{\mathbb{E}v_\infty\left(\tilde{W}_0, \tilde{X}_0, \tau - \tau_0\right)}* 100.
    \label{eq:7.33}
\end{equation}
Contour plots are used to explore $D$'s dependence as a function of selected pairs of parameters for the system. In the top right plot of Figure~\ref{fig:7.1}, the improvement increases as $\gamma$ grows closer to 1. As $\gamma$ grows, previous terms in the $\epsilon_n(k)$ sequence play a larger role so the value of forecasts increases. On the other hand, as $g$ grows closer to 1, the dependence across time of the $W_n$ sequence grows, and so the value of forecasts decreases. On the top right plot, the value of forecasts increases as the dependence of $X_{n+1}$ on $X_n$ grows (with $\rho$). In the bottom left plot, the value of forecasts increases with the noise in the weather sequence $W_n$ (controlled by $\sigma^2$) but decreases with the noise in the building temperature sequence $X_n$ (controlled by $\sigma_V^2$). Finally, the bottom right plot shows a symmetry in the improvement with respect to the control setpoint $\tau$ around $\tau_0$. We note that $\tau_0=\mathbb{E}X_0$, the building temperature in the absence of control. The value of forecasts decreases when the target setpoint is farther from $\tau_0$: as the weight of the action in the value function grows, the relative value of forecasts weakens.

\section{Conclusion}
\label{sec:8}
In this work we introduced the first principled and mathematically consistent framework for the incorporation of forecasts into MDPs in the setting of state space models with linear dynamics, using no ad hoc elements to add forecast information into the MDP setting. In this framework, we discussed the different ways in which forecast information can be incorporated (static, dynamic, static and dynamic together). Through an illustrative energy system control example, we provided a numerical comparison of the optimal value functions for the setting with no forecasts to the setting with dynamic forecasts.
\par The introduction of this framework opens the door to several theoretical and applied research questions, \textit{e.g.} on how the quality of forecasts affects control methods in different disciplines and in different applications. Potential theoretical research directions include extensions to periodic Markov chains (e.g. to model time-of-day effects), non-stationary Markov chains, forecast
updates that are not synchronized with decisions epochs, and Markov chains with nonlinear dynamics.

\vskip6pt

\enlargethispage{20pt}

\aucontribute{JAC and PWG conceived of the study and wrote the manuscript.}

\competing{The authors declare that they have no competing interests.}

\funding{Funding for this research was supported by a State Grid Graduate Student Fellowship through the Stanford Bits \& Watts initiative and by Total SE.}

\ack{The authors thank the members of the Stanford Cooler project, in particular Professor Sally Benson and Rui Yan, for valuable insights and discussions, as well as the two anonymous referees whose comments helped us to significantly improve our work.}

\bibliographystyle{rsta}
\bibliography{biblio}

\end{document}